\documentclass[12pt]{amsart}
\usepackage{graphicx, amssymb}
\numberwithin{equation}{section}

\newtheorem{Thm}{Theorem}[section]

\newtheorem{Lem}[Thm]{Lemma}

\def\bN {\mathbf{N}}
\def\bR {\mathbf{R}}
\def\bS {\mathbf{S}}
\def\bT {\mathbf{T}}
\def\bZ {\mathbf{Z}}

\def\cC {\mathcal{C}}

\def\cK {\mathcal{K}}

\def\cR {\mathcal{R}}

\def\g {{\gamma}}

\def\eps {{\epsilon}}

\def\si {{\sigma}}
\def\Ups {{\Upsilon}}
\def\om {{\omega}}
\def\Om {{\Omega}}

\def\rstr {{\big |}}
\def\indc {{\bf 1}}

\def\la {\langle}
\def\ra {\rangle}

\def\d {{\partial}}
\def\grad {{\nabla}}

\newcommand{\Supp}{\operatorname{supp}}

\newcommand{\Dist}{\operatorname{dist}}

\begin{document}
\title[Periodic Lorentz Gas]
      {On the Periodic Lorentz Gas\\ and the Lorentz Kinetic Equation}

\date{\today}

\author[F. Golse]{Fran\c cois Golse}
\address[F. G.]%
{Ecole polytechnique, CMLS, F91128 Palaiseau cedex 
\& Laboratoire Jacques-Louis Lions 
 Bo\^\i te courrier 187, F75252 Paris cedex 05}
\email{golse@math.polytechnique.fr}

\begin{abstract}
We prove that the Boltzmann-Grad limit of the Lo- rentz gas with periodic distribution of scatterers
cannot be described with a linear Boltzmann equation. This is at variance with the case of a Poisson
distribution of scatterers, for which the convergence to the linear Boltzmann equation was proved
by Gallavotti [{\it Phys. Rev. (2)} {\bf 185}, 308 (1969)]. The arguments presented here complete
the analysis in [Golse-Wennberg, {\it M2AN Mod\'el. Math. et Anal. Num\'er.} {\bf 34}, 1151 (2000)],
where the impossibility of a kinetic description was established only in the case of absorbing
obstacles. The proof is based on estimates on the distribution of free-path lengths established in
[Golse-Wennberg {\it loc.cit.}] and in [Bourgain-Golse-Wennberg, {\it Commun. Math. Phys.} {\bf 190},
491 (1998)], and on a classical result on the spectrum of the linear Boltzmann equation which can be
found in [Ukai-Point-Ghidouche, {\it J. Math. Pures Appl. (9)} {\bf 57}, 203 (1978)].
\end{abstract}

\maketitle

\noindent
\textbf{Key-words:} Lorentz gas; Sinai billiards; Kinetic theory; Boltzmann-Grad limit.

\noindent
\textbf{MSC:} 82C40 (37A60,37D50) 

\section{Introduction}

About 100 years ago, Lorentz \cite{Lo} proposed the following linear kinetic equation to describe
the motion of electrons in a metal:
\begin{equation}
\label{KntLrtz}
(\d_t+v\cdot\grad_x+\tfrac1m F(t,x)\cdot\grad_v)f(t,x,v)=N_{at}r_{at}^2|v|\cC(f(t,x,\cdot))(v)
\end{equation}
where $f(t,x,v)$ is the (phase space) density of electrons which, at time $t$, are located at $x$
and have velocity $v$. In Eq.~(\ref{KntLrtz}), $F$ is the electric force field, $m$ the mass of
the electron, while $N_{at}$ and $r_{at}$ designate respectively the number of metallic atoms per
unit volume and the radius of each such atom. Finally $\cC(f)$ is the collision integral: it acts
on the velocity variable only, and is given, for all continuous $\phi\equiv\phi(v)$ by the formula
\begin{equation}
\label{CllsLrntz}
\cC(\phi)(v)=\int_{|\om|=1,v\cdot\om>0}\bigl(\phi(v-2(v\cdot\om)\om)-\phi(v)\bigr)\cos(v,\om)d\om\,.
\end{equation}
In the case where $F\equiv 0$, Gallavotti \cite{Ga1,Ga2} proved that Eq.~(\ref{KntLrtz}) describes
the Boltzmann-Grad limit of a gas of point particles undergoing elastic collisions on a random
(Poisson) configuration of spherical obstacles. His result was successively strengthened by Spohn
\cite{Sp}, and by Boldrighini-Bunimovich-Sinai \cite{BBS}.

In the presence of an external, non-zero electric force $F$ and for the same random configuration
of absorbing obstacles as in \cite{Ga1,Ga2}, Desvillettes-Ricci \cite{DeRi} proved recently that
the Boltzmann-Grad limit of a gas of point particles leads to a non-Markovian equation --- see
also an earlier, similar observation by Bobylev-Hansen-Piasecki-Hauge \cite{BHPH}.

The case of periodic configuration of obstacles, perhaps closer to Lorentz's original ideas, completely 
differs from the random case. In the case of absorbing obstacles, and without external force $F$,
several results suggest that the Boltzmann-Grad limit is non-Markovian \cite{BGW,GW,CG}. However,
all these results are based on explicit computations that are possible only in the case of absorbing
obstacles.

In the present note, we show that neither equation (\ref{KntLrtz}) nor any variant thereof can describe
the Boltzmann-Grad limit of the periodic Lorentz gas with no external force ($F\equiv 0$) and in the
case of reflecting obstacles.

The crucial observation (already made in \cite{BGW,GW}) is that the distribution of first hitting times
satisfies, in the periodic case, an inequality stated below as Theorem \ref{TH-BGW} that prevents these
first hitting times from begin exponentially distributed. The probabilistic representation of equation
(\ref{KntLrtz}) by a jump process with exponentially distributed jump times plus a drift (see for instance
\cite{Pa}) suggests that the Boltzmann-Grad limit of the periodic Lorentz gas with reflecting obstacles
cannot be described by (\ref{KntLrtz}).

In the present paper, we give a complete proof of this fact (see Theorem \ref{TH-MN} below), which by 
the way does not appeal to the probabilistic representation of (\ref{KntLrtz}). In addition, the method of 
proof used here provides an explicit estimate of the difference between the single particle phase space 
density for the periodic Lorentz gas in the Boltzmann-Grad limit and the solution of the Lorentz kinetic 
equation (\ref{KntLrtz}) --- see the inequality (\ref{Inq-r}) and the discussion thereafter. This result has
been announced in \cite{GolICMP}

Recently, Ricci and Wennberg \cite{RiWe} have considered the following interesting variant of the
periodic Lorentz gas studied here. Their microscopic model consists of gas of point particles in a 
periodic configuration of obstacles that are bigger than in the Boltzmann-Grad scaling considered 
in the present paper. However, some of these obstacles are removed with a probability carefully
chosen in terms of the obstacle radius so that the mean collision time remains of order one. They
proved that the expected single-particle phase-space density so obtained converges to a solution
of the Lorentz kinetic equation (\ref{KntLrtz}). This result suggests that the failure of (\ref{KntLrtz})
to capture the Boltzmann-Grad limit of the periodic Lorentz gas is a very unstable phenomenon 
that is specific to the periodic setting and likely to disappear whenever some amount of randomness
is injected in the microscopic system.

\section{The periodic Lorentz gas}

Let $D\in\bN$, $D\ge 2$. For each $r\in(0,\tfrac12)$, we consider the domain
$$
Z_r=\{x\in\bR^D\,|\,\Dist(x,\bZ^D)>r\}\,,
$$
which is usually referred to as ``the billiard table".

The {\it Lorentz gas} is the dynamical system corresponding to a cloud of point particles that move
freely in $Z_r$ --- collisions between particles being neglected --- and are specularly reflected at
the boundary of each obstacle --- the obstacles being the balls of radius $r$ centered at the lattice
points, i.e. the connected components of $Z_r^c$. Since collisions between particles are neglected,
one can equivalently consider a single particle whose initial position and velocity are appropriately
distributed in $Y_r\times\bS^{D-1}$.

In this dynamical system, a particularly important notion is that of ``free path length": see the
review article by Bunimovich \cite{Bu}, pp. 221-222.

The free path length --- or ``(forward) exit time" --- for a particle starting from $x\in Z_r$ in
the direction $v\in\bS^{D-1}$ is defined as
\begin{equation}
\label{DfntTau}
\tau_r(x,v)=\inf\{t>0\,|\,x+tv\in\d Z_r\}\,.
\end{equation}
The function $\tau_r$ is then extended by continuity to the non-characteristic part of the boundary
of the phase-space, i.e. to
$$
\{(x,v)\in\d Z_r\times\bS^{D-1}\,|\,v\cdot n_x\not=0\}\,,
$$
where $n_x$ designates the inward unit normal field on $\d Z_r$ (pointing toward $Z_r$, i.e. outside
of the scatterers). Because $Z_r$ is invariant under $\bZ^D$-translations, one has
$$
\tau_r(x+k,v)=\tau_r(x,v)\hbox{ for each }(x,v)\in Z_r\times\bS^{D-1}\hbox{ and }k\in \bZ^D\,.
$$
Hence $\tau_r$ can be seen as a $[0,+\infty]$-valued function defined on $Y_r\times\bS^{D-1}$ (and
a.e. on $\overline{Y}_r\times\bS^{D-1}$), where $Y_r=Z_r/\bZ^D$.

Whenever the components of $v\in\bS^{D-1}$ are rationally independent --- i.e. if $k\cdot v\not=0$
for each $k\in\bZ^D\setminus\{0\}$ --- each orbit of the linear flow $x\mapsto x+tv$ is dense on
$\bT^D=\bR^D/\bZ^D$, and thus $\tau_r(x,v)<+\infty$ for each $x\in Z_r$.

\begin{figure}
\centering
\includegraphics[height=5in]{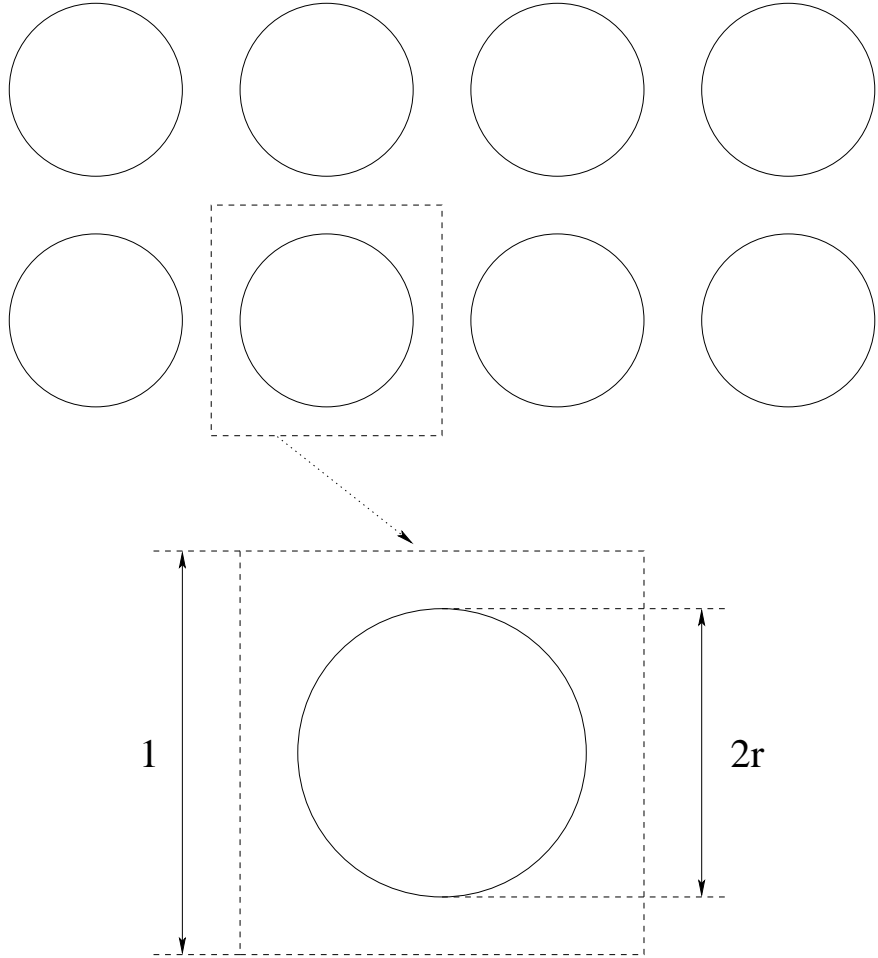}
\caption{The billiard table $Z_r$ and the punctured torus $Y_r$ \label{fig1}}
\end{figure}

On the measurable space $Y_r\times\bS^{D-1}$ equipped with its Borel $\si$-algebra, we define $\mu_r$
as the probability measure proportional to the Lebesgue measure on $Y_r \times\bS^{D-1}$: in other
words\footnote{If $A$ is a measurable $d$-dimensional set in $\bR^D$ ($d\le D$), we denote by $|A|$
its $d$-dimensional volume. Here $|Y_r|$ is the Lebesgue measure of any fundamental domain of the
quotient space $Y_r$, i.e. of the unit cube with a ball of radius $r$ removed.}
$$
d\mu_r(y,v)=\frac{dydv}{|Y_r||\bS^{D-1}|}\,.
$$
Define the distribution of $\tau_r$ under $\mu_r$ by the usual formula
$$
\Phi_r(t):=\mu_r\left(\{(y,v)\in Y_r\times\bS^{D-1}\,|\,\tau_r(y,v)\ge t\right\})\,.
$$

\begin{Thm}\label{TH-BGW}
Let $D\ge 2$. There exist two positive constants $C_1$ and $C_2$ such that, for each $r\in(0,\tfrac12)$
and each $t>1/r^{D-1}$
$$
\frac{C_1}{tr^{D-1}}\le\Phi_r(t)\le\frac{C_2}{tr^{D-1}}\,.
$$
\end{Thm}

In the theorem above, the lower bound in the case $D=2$ and the upper bound for all $D\ge 2$ were
proved by Bourgain-Golse-Wennberg \cite{BGW}; the lower bound was extended to the case of any $D\ge 2$
by Golse-Wennberg \cite{GW}. More precise results concerning $\Phi_r(t/r)$ in space dimension $D=2$
have recently been obtained by Caglioti-Golse \cite{CG} and Boca-Zaharescu \cite{BZ}; however, only the
result above (Theorem \ref{TH-BGW}) is used in the present paper.

\section{The linear Boltzmann equation}

The linear Boltzmann equation for a free\footnote{I.e. without external force.} gas of particles moving
at speed $1$ is
\begin{equation}
\label{Bltz}
(\d_t+v\cdot\grad_x)f(t,x,v)=\cC(f(t,x,\cdot))(v)
\end{equation}
where $f\equiv f(t,x,v)$ is the single particle phase-space density, also known as ``distribution function". 
In other words, $f(t,x,v)$ is the density of particles which, at time $t$, are located at $x\in\bR^D$ and 
move in the direction $v\in\bS^{D-1}$. The term $\cC(f)$ is the collision integral; specifically, $\cC$ is
an integral operator of the form
\begin{equation}
\label{CllsBltz}
\cC(\phi)(v)=\si\int_{\bS^{D-1}} k(v,w)(\phi(w)-\phi(v))dw
\end{equation}
where $\si>0$ and $k\in C(\bS^{D-1}\times\bS^{D-1})$ satisfies the following properties
\begin{equation}
\label{Prpt-k}
k(v,w)=k(w,v)>0\,,\quad\int_{\bS^{D-1}}k(v,w)dv=1\,,\hbox{ for all }v,w\in\bS^{D-1}\,.
\end{equation}
Changing variables according to $\om\mapsto w=v-2(v\cdot\om)\om$ in the Lorentz collision integral for
$D=3$ shows that
$$
\begin{aligned}
\tfrac1\pi\int_{\om\in\bS^2,v\cdot\om>0}&\bigl(\phi(v-2(v\cdot\om)\om)-\phi(v)\bigr)\cos(v,\om)d\om
\\
&=
\int_{\bS^{2}}k(v,w)(\phi(w)-\phi(v))dw\hbox{ with }k(v,w)\equiv\tfrac1{4\pi}\,.
\end{aligned}
$$
Hence, the Lorentz kinetic model (\ref{KntLrtz})-(\ref{CllsLrntz}) is a special case of the linear
Boltzmann equation (\ref{Bltz})-(\ref{CllsBltz}).

In the sequel, we shall restrict our attention to the case where the linear Boltzmann equation is
posed in a periodic box. Without loss of generality, we assume that this periodic box has size 1,
and that $\bS^{D-1}$ is endowed with its rotationally invariant unit measure, henceforth denoted 
by $dv$. Finally we denote by $\la\cdot\ra$ the average with respect to both variables $x$ and 
$v$:
\begin{equation}
\label{Dfnt<>}
\la\phi\ra=\iint_{\bT^D\times\bS^{D-1}}\phi(x,v)dxdv\,.
\end{equation}

Consider then the unbounded operator $A$ on $L^2(\bT^D\times\bS^{D-1})$ defined by
$$
(A\phi)(x,v)=-v\cdot\grad_x\phi(x,v)-\si\phi(x,v)+\si\int_{\bS^{D-1}}k(v,w)\phi(x,w)dw
$$
with domain
$$
D(A)=\{\phi\in L^2(\bT^D\times\bS^{D-1})\,|\,v\cdot\grad_x\phi\in L^2(\bT^D\times\bS^{D-1})\}\,.
$$

We recall the following result, originally proved in \cite{UPG} for the more complicated case
of the linearization of Boltzmann's equation at a uniform Maxwellian state.

\begin{Thm}\label{TH-SP}
Let $\si>0$ and $k\in C(\bS^{D-1}\times\bS^{D-1})$ satisfy the assumptions (\ref{Prpt-k}). Then the
operator $A$ generates a strongly continuous contraction semigroup on $L^2(\bR^D\times\bS^{D-1})$,
and there exists positive constants $c$ and $\g$ such that
$$
\|e^{tA}\phi-\la\phi\ra\|_{L^2(\bT^D\times\bS^{D-1})}
    \le ce^{-\g t}\|\phi\|_{L^2(\bT^D\times\bS^{D-1})}\,,\quad t\ge 0\,,
$$
for each $\phi\in L^2(\bT^D\times\bS^{D-1})$.
\end{Thm}

\section{The non-convergence result}

Throughout this section, we denote by $\eps$ the sequence $(1/n)_{n\in\bN^*}$. For some given
$r_*\in(0,\tfrac12)$, we set $r=r_*\eps^{1/(D-1)}$, and we denote by $\Om_\eps$ the open subset
$\eps Z_r$ with this particular choice of $r$.

\subsection{Reflecting vs. absorbing obstacles}

First we define the billiard flow on the scaled billard table $\Om_\eps$. It is a one-parameter
group on $\Om_\eps\times\bS^{D-1}$ denoted by
\begin{equation}
\label{DfntS}
S^\eps_t:\,
\Om_\eps\times\bS^{D-1}\ni(x,v)\mapsto (X^\eps_t(x,v),V^\eps_t(x,v))\in\Om_\eps\times\bS^{D-1}
\end{equation}
and defined in the following manner:

\begin{itemize}
\item $S^\eps_0(x,v)=(x,v)$,

\item if $t$ is such that $X^\eps_t(x,v)\notin\d\Om_\eps$, then $S^\eps_t(x,v)$ is smooth in the
time variable near that particular value of $t$ and one has
\begin{equation}
\label{DfntS1}
\frac{d}{dt}{X}^\eps_t(x,v)=V^\eps_t(x,v)\,,\quad\frac{d}{dt}{V}^\eps_t(x,v)=0\,;
\end{equation}

\item if $t$ is such that $X^\eps_t(x,v)=x^*\in\d\Om_\eps$, then $S^\eps_t(x,v)$ has the following
jump discontinuity (in the velocity component only)
\begin{equation}
\label{DfntS2}
X^\eps_{t+0}(x,v)=X^\eps_{t-0}(x,v)\,,\quad V^\eps_{t+0}(x,v)=\cR(n_{x^*})V^\eps_{t-0}(x,v)\,,
\end{equation}
where $\cR(n_{x^*})$ is the symmetry with respect to the hyperplane orthogonal to the inner unit
normal $n_{x^*}$ at $x^*\in\d\Om_\eps$: in other words
\begin{equation}
\label{DfntRn}
\cR(n_{x^*})\xi=\xi-2(\xi\cdot n_{x^*})n_{x^*}\,,\quad\xi\in\bR^D\,.
\end{equation}
\end{itemize}

\begin{figure}
\centering
\includegraphics[height=3.5in]{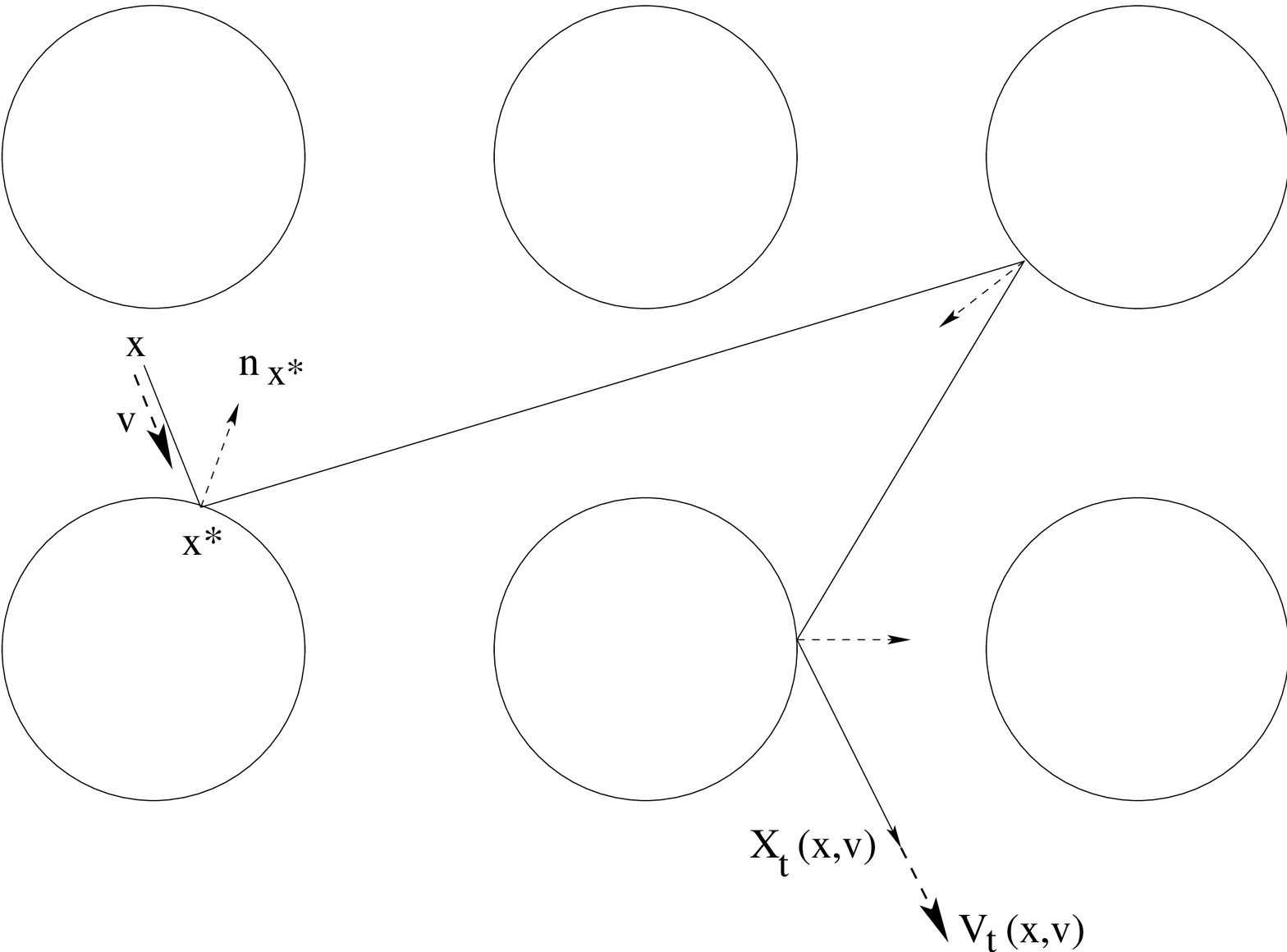}
\caption{The billiard dynamics on $\Om_\eps$\label{fig2}}
\end{figure}

Notice that this dynamics is mechanically reversible, i.e. for all $t\in\bR$, one has
$$
S^\eps_{-t}(x,v)=S^\eps_t(x,-v)\,,\quad(x,v)\in\Om_\eps\times\bS^{D-1}\,.
$$
For each $t\in\bR$, the map $S^\eps_t$ is invariant under $\bZ^D$-translations, in the following
sense: for all $k\in\bZ^D$ and $(x,v)\in\Om_\eps\times\bS^{D-1}$, one has
$$
S^\eps_t(x+k,v)=(X^\eps_t(x,v)+k,V^\eps_t(x,v))=S^\eps_t(x,v)+(k,0)\,.
$$
Hence $S^\eps_t$ defines a one-parameter group --- still denoted by $S^\eps_t$ --- on the quotient
space $\Ups_\eps\times\bS^{D-1}$, where $\Ups_\eps=\Om_\eps/\bZ^D$, and the restriction of the 
measure $dxdv$ (defined on $\bT^D\times\bS^{D-1}$) to $\Ups_\eps\times\bS^{D-1}$ is invariant 
under $S^\eps_t$.

Starting from the billiard flow $S^\eps_t$ on $\Ups_\eps\times\bS^{D-1}$, we define a unitary group
$\hat S^\eps_t$ on $L^2(\Ups_\eps\times\bS^{D-1},dxdv)$ by the formula
\begin{equation}
\label{Dfnt-hS}
\hat S^\eps_tf(x,v)=f(S^\eps_{-t}(x,v))\,,\quad(x,v)\in\Ups_\eps\times\bS^{D-1}\,,\,\,t\in\bR\,.
\end{equation}

Next we turn to the case of the same billiard table, but with absorbing obstacles. By this, we mean
that whenever a particle hits the boundary of $\Om_\eps$, it disappears. Equivalently, one may think
of $\Om_\eps$ as a sieve, with particles falling into the holes (the components of $\Om_\eps^c$).

The analogue of $\hat S^\eps_t$ in the case of absorbing obstacles is the contraction semigroup
$T^\eps_t$ defined on $L^2(\Ups_\eps\times\bS^{D-1},dxdv)$ by the formula
\begin{equation}
\label{DfntT}
T^\eps_tg(x,v)=g(x-tv,v)\indc_{\tau_r(x/\eps,-v)>t/\eps}\,,
    \quad(x,v)\in\Ups_\eps\times\bS^{D-1}\,,\,\,t\in\bR\,,
\end{equation}
where $\tau_r$ is the free path length defined in (\ref{DfntTau}) --- recall that $r_*\eps^{1/(D-1)}$.

Finally, $L^2(\Ups_\eps\times\bS^{D-1},dxdv)$ can be embedded as a subspace of 
$L^2(\bT^D\times\bS^{D-1},dxdv)$, identifying each function defined a.e. on $\Ups_\eps\times\bS^{D-1}$
with its extension by $0$ in the complement of $\Ups_\eps\times\bS^{D-1}$.

The maps $\hat S^\eps_t$ and $T^\eps_t$ are related by the following elementary inequality: for each
$f\in L^2(\Ups_\eps\times\bS^{D-1},dxdv)$ and each $t\ge 0$,
\begin{equation}
\label{InqST}
f\ge 0\hbox{ a.e. implies that }\hat S^\eps_tf(x,v)\ge T^\eps_tf(x,v)\hbox{ a.e. }
\end{equation}
on $\Ups_\eps\times\bS^{D-1}$. Indeed, the formulas (\ref{DfntT}), (\ref{Dfnt-hS}), and (\ref{DfntS1})
imply that
$$
\begin{aligned}
T^\eps_tf(x,v)=f(x-tv,v)\indc_{\tau_r(x/\eps,-v)>t/\eps}
=f(S^\eps_{-t}(x,v))\indc_{\tau_r(x/\eps,-v)>t/\eps}
\\
\le f(S^\eps_{-t}(x,v))=\hat S^\eps_tf(x,v)\,.
\end{aligned}
$$
Finally, we give the PDE interpretation of the operators $\hat S^\eps_t$ and $T^\eps_t$. The
function $F_\eps(t,x,v)=\hat S^\eps_tf(x,v)$ is the solution of
$$
\begin{aligned}
\d_tF_\eps+v\cdot\nabla_xF_\eps&=0\,,&&\quad (x,v)\in\Ups_\eps\times\bS^{D-1}\,,
\\
F_\eps(t,x,v)&=F_\eps(t,x,\cR(n_x)v)\,,&&\quad (x,v)\in\d\Ups_\eps\times\bS^{D-1}\,,
\\
F_\eps\rstr_{t=0}&=f\,,&&
\end{aligned}
$$
(where $\cR(n_x)$ designates the reflection with respect to the hyperplane orthogonal to $n_x$: see fla.
(\ref{DfntRn})), while the function $G_\eps(t,x,v)=\hat T^\eps_tg(x,v)$ is the solution of
$$
\begin{aligned}
\d_tG_\eps+v\cdot\nabla_xG_\eps&=0\,,&&\quad (x,v)\in\Ups_\eps\times\bS^{D-1}\,,
\\
G_\eps(t,x,v)&=0\,,&&\quad x\in\d\Ups_\eps\,,\,\,v\cdot n_x>0\,,
\\
G_\eps\rstr_{t=0}&=g\,.&&
\end{aligned}
$$

\subsection{Main result}

The paper by Lorentz \cite{Lo} described in the introduction suggests the following question, in
the case of space dimension $D=3$:

``Let $f^{in}\in L^2(\bT^3\times\bS^2)$. Does $\hat S^\eps_t(f^{in}\indc_{\Ups_\eps\times\bS^2})$
(or any subsequence thereof) converge in $L^\infty((0,+\infty)\times\bT^3\times\bS^2)$ weak-* as
$\eps\to 0$ to the solution $f$ of (\ref{KntLrtz}) on $\bT^3\times\bS^2$ with $N_{at}=1$, $r_{at}
=r_*$, and with initial data $f^{in}$?"

The answer to that question is negative, as shown by the following theorem.

\begin{Thm}\label{TH-MN}
Assume that the space dimension satisfies $D\ge 2$. There exist initial data 
$f^{in}\in L^2(\bT^D\times\bS^{D-1})$ such that, for any $\si>0$ and any function 
$k\in C(\bS^{D-1}\times\bS^{D-1})$ satisfying (\ref{Prpt-k}), no subsequence of 
$\hat S^\eps_t(f^{in}\indc_{\Ups_\eps\times\bS^{D-1}})$ converges in 
$L^\infty((0,+\infty)\times\bT^D\times\bS^{D-1})$ weak-* to the solution of
\begin{equation}
\label{LnBltz}
(\d_t+v\cdot\grad_x)f(t,x,v)=\si\int_{\bS^{D-1}}k(v,w)\left(f(t,x,w)-f(t,x,v)\right)dw
\end{equation}
on $\bT^D\times\bS^{D-1}$, with initial data
\begin{equation}
\label{LnBltzCndn}
f\rstr_{t=0}=f^{in}\,.
\end{equation}
\end{Thm}

\subsection{Proof of Theorem \ref{TH-MN}.}

We define
$$
\Phi_\eps(t,x,v)=\indc_{\eps\tau_r(x/\eps,-v)>t}\,.
$$
Since $0\le\Phi_\eps\le 1$, the sequence $\Phi_\eps$ is relatively weakly-* compact in 
$L^\infty(\bR_+\times\bT^D\times\bS^{D-1})$. From now on, we consider a subsequence 
of $\eps=1/n$, denoted $\eps'$, so that $\Phi_{\eps'}\to\Phi$ in 
$L^\infty(\bR_+\times\bT^D\times\bS^{D-1})$ weak-*, and denote $r'=r_*{\eps'}^{1/{D-1}}$. 
Then, the function $\Phi$ is independent of $x$, i.e.
\begin{equation}
\label{LmtPhi}
\Phi\equiv\Phi(t,v)\,,
\end{equation}
as shown by the following classical lemma.

\begin{Lem} \label{LM-LMPH}
Let $Z$ be a separable locally compact metric space, endowed with a Borel measure 
$m$. Let $u_n(x,z)\equiv U_n(nx,z)$, where $(U_n)_{n\ge 1}$ is a bounded sequence 
of elements of $L^\infty(\bT^D\times Z)$. Any weak-* limit point of the sequence 
$(u_n)_{n\ge 1}$ as $n\to+\infty$ is a function of $z$ alone (i.e. independent of $x$).
\end{Lem}

We postpone the proof of Lemma \ref{LM-LMPH} until the end of this section.

In the sequel, we consider an arbitrary initial data $\rho\equiv\rho(x)\in L^\infty(\bT^D)$
that is independent of $v$ and a.e. nonnegative. We define
\begin{equation}
\label{Dfnt-fe}
f_\eps(t,x,v)=\hat S^\eps_t(\rho\indc_{\Ups_\eps\times\bS^{D-1}})
    =\rho(X^\eps_{-t}(x,v))\indc_{\Ups_\eps}(x)
\end{equation}
and
\begin{equation}
\label{Dfnt-ge}
g_\eps(t,x,v)=T^\eps_t(\rho\indc_{\Ups_\eps\times\bS^{D-1}})=
    \rho(x-tv)\Phi_\eps(t,x,v)\indc_{\Ups_\eps}(x)\,.
\end{equation}
Since $\indc_{\Ups_\eps}\to 1$ a.e. on $\bT^D\times\bS^{D-1}$ and $|\indc_{\Ups_\eps}|\le 1$,
one has, by dominated convergence,
\begin{equation}
\label{lmt-ge}
g_{\eps'}\to\rho(x-tv)\Phi(t,v)
    \hbox{ in }L^\infty(\bR_+\times\bT^D\times\bS^{D-1})\hbox{ weak-*}\,.
\end{equation}
On the other hand, (\ref{Dfnt-fe}) shows that $\|f_\eps\|_{L^\infty}\le\|\rho\|_{L^\infty}$;
therefore, possibly after extraction of a subsequence (still denoted $\eps'$), one has $f_{\eps'}
\to f$ in $L^\infty(\bR_+\times\bT^D\times\bS^{D-1})$ weak-*. We recall that $\rho\ge 0$ a.e. on
$\bT^D$. Therefore, because of the inequality (\ref{InqST}) and of the weak-* limit (\ref{lmt-ge}),
one has
\begin{equation}
\label{lmt-fe}
f(t,x,v)\ge\rho(x-tv)\Phi(t,v)\,.
\end{equation}
In particular,
\begin{equation}
\label{LwrBnd}
\begin{aligned}
\iint_{\bT^D\times\bS^{D-1}}f(t,x,v)^2dxdv
&\ge
\iint_{\bT^D\times\bS^{D-1}}\rho(x-tv)^2\Phi(t,v)^2dxdv
\\
&=
\int_{\bT^D}\rho(y)^2dy\int_{\bS^{D-1}}\Phi(t,v)^2dv
\\
&\ge
\|\rho\|^2_{L^2}\left(\int_{\bS^{D-1}}\Phi(t,v)dv\right)^2\,,
\end{aligned}
\end{equation}
where the last inequality follows from Jensen's inequality. By assumption, $\Phi_{\eps'}\to\Phi$
in $L^\infty(\bR_+\times\bT^D\times\bS^{D-1})$ weak-*; because of (\ref{LmtPhi}), one has
$$
\iint_{\bT^D\times\bS^{D-1}}\Phi_{\eps'}(t,x,v)dxdv\to\int_{\bS^{D-1}}\Phi(t,v)dv
\hbox{ in }L^\infty(\bR_+)\hbox{ weak-*.}
$$
Furthermore, by Theorem \ref{TH-BGW},
$$
\begin{aligned}
\iint_{\bT^D\times\bS^{D-1}}&\Phi_{\eps'}(t,x,v)dxdv
\\
&=
\hbox{meas}(\{(x,v)\in\Ups_{\eps'}\times\bS^{D-1}\,|\,\tau_{r'}(x/\eps',-v)>t/\eps'\})
\\
&=
\hbox{meas}(\{(y,v)\in Y_{r'}\times\bS^{D-1}\,|\,\tau_{r'}(y,-v)>t/\eps'\})
\\
&\ge
\frac{C_1}{\frac{t}{\eps'}{r'}^{D-1}}=\frac{C_1}{tr_*^{D-1}}\,,\quad\hbox{ for all }t>1/r_*^{D-1}\,.
\end{aligned}
$$
Hence, the inequality (\ref{LwrBnd}) becomes
\begin{equation}
\label{iintf>}
\left(\iint_{\bT^D\times\bS^{D-1}}f(t,x,v)^2dxdv\right)^{1/2}\ge\frac{C_1}{tr_*^{D-1}}\|\rho\|_{L^2}\,,
\quad\hbox{ for all }t>1/r_*^{D-1}\,.
\end{equation}

Assume that $f$ is the solution to (\ref{LnBltz})-(\ref{LnBltzCndn}) with $f^{in}=\rho$. 
By Theorem \ref{TH-SP}, one has
$$
\left\|f(t,\cdot,\cdot)-\int_{\bT^D}\rho(x)dx\right\|_{L^2(\bT^D\times\bS^{D-1})}
\le
ce^{-\g t}\|\rho\|_{L^2}
$$
for all $t\ge 0$. In particular, for all $t\ge 0$, one has
\begin{equation}
\label{iintf<}
\left(\iint_{\bT^D\times\bS^{D-1}}f(t,x,v)^2dxdv\right)^{1/2}
\le
\int_{\bT^D}\rho(x)dx+ce^{-\g t}\|\rho\|_{L^2}\,.
\end{equation}

In conclusion, if $f$ is the solution to (\ref{LnBltz})-(\ref{LnBltzCndn}), then the initial data
$\rho$ (assuming it is not a.e. $0$) must satisfy the inequality
\begin{equation}
\label{Inq-r}
\frac{C_1}{tr_*^{D-1}}\le\frac{\|\rho\|_{L^1}}{\|\rho\|_{L^2}}+ce^{-\g t}\,,
\quad\hbox{ for all }t>1/r_*^{D-1}\,.
\end{equation}
At this point, we recall that this inequality holds for any arbitrary $\rho\in L^\infty(\bT^D)$ such
that $\rho\ge 0$ a.e.: hence the ratio $\|\rho\|_{L^1}/\|\rho\|_{L^2}$ can be made arbitrarily small.
For instance, one can choose $\rho$ as follows: pick $b$, a bump function on $\bR^D$ satisfying
$$
0\le b\le 1\,,\quad\Supp(b)\subset[-\tfrac14,\tfrac14]^D\,.
$$
For $m\in\bN^*$, define $\rho$ to be the unique $\bZ^D$-periodic function such that
$$
\rho\rstr_{[-\tfrac12,\tfrac12)^D}:\,x\mapsto b(mx)\,.
$$
Then
$$
\|\rho\|_{L^1(\bT^D)}=m^{-D}\|b\|_{L^1(\bR^D)}\,,\quad \|\rho\|_{L^2(\bT^D)}=m^{-D/2}\|b\|_{L^2(\bR^D)}\,,
$$
so that
$$
\frac{\|\rho\|_{L^1(\bT^D)}}{\|\rho\|_{L^2(\bT^D)}}=m^{-D/2}\to 0\hbox{ as }m\to+\infty\,.
$$
Since the ratio $\|\rho\|_{L^1}/\|\rho\|_{L^2}$ can be made arbitrarily small, the inequality
(\ref{Inq-r}) would entail
$$
\frac{C_1}{tr_*^{D-1}}\le ce^{-\g t}\,,\quad\hbox{ for all }t>1/r_*^{D-1}\,,
$$
which is manifestly wrong for $t$ large enough --- specifically, larger than the unique zero of the
function $t\mapsto C_1e^{\g t}-cr_*^{D-1}t$, denoted by $t_*(C_1,c,\g,r_*)$.

Hence, the assumption that, for each nonnegative $\rho\equiv\rho(x)\in L^\infty(\bT^D)$, there exists
a subsequence of $\hat S^\eps_t(\rho\indc_{\Ups_\eps\times\bS^{D-1}})$ converging in $L^\infty$ 
weak-* to the solution of (\ref{LnBltz})-(\ref{LnBltzCndn}) with $f^{in}=\rho$ is wrong. This concludes 
the proof of Theorem \ref{TH-MN}, once Lemma \ref{LM-LMPH} is proved.

\subsection{Proof of Lemma \ref{LM-LMPH}}

Assume that
$$
u_{n_q}\to u\hbox{ in }L^\infty(\bT^D\times Z)\hbox{ weak-* as }n_q\to+\infty\,.
$$
For each $k\in\bZ^D$, one has
$$
\begin{aligned}
\hat u_{n_q}(k,z)=\int_{\bT^D}e^{-i2\pi k\cdot x}u_{n_k}(x,z)dx
&=\hat U_{n_q}(k,z)&&\hbox{ if }n_q|k\,,
\\
&=0&&\hbox{ otherwise.}
\end{aligned}
$$
Since the assumed convergence entails
$$
\hat u_{n_q}(k,\cdot)\to\hat u(k,\cdot)\hbox{ in }L^\infty(Z)\hbox{ weak-* for each }k\in\bZ^D\,,
$$
this shows that
$$
\hat u(k,\cdot)=0\hbox{ unless }k=0\,.
$$
Hence $u$ is independent of $x$, as announced.

\section{Final remarks}

We have demonstrated the impossibility of representing the Boltz- mann-Grad limit of the 
periodic Lorentz gas by a linear Boltzmann equation.

Notice that this impossibility results solely from the lower bound on the distribution of free
path lengths in Theorem \ref{TH-BGW}. This is in fact not too surprising since the 
probabilistic representation of equation (\ref{KntLrtz}) involves in particular exponentially 
distributed jump times for the velocity process.

Notice also the choice of initial data in Theorem \ref{TH-MN}. Since one knows that 
particles moving in appropriately chosen rational directions may not encounter any
obstacles on the billiard table $Z_r$, it may seem somewhat surprising that the 
contradiction in Theorem \ref{TH-MN} is obtained by considering isotropic initial data 
(i.e. initial data that are independent of the angle variable) instead of pencils of 
particles concentrated in phase-space on those rational directions that avoid all 
obstacles. In fact, the contribution of such rational directions is already taken into 
account in the lower bound in Theorem \ref{TH-BGW}. Besides, the case of isotropic 
initial densities is somewhat more natural in the context of Theorem \ref{TH-MN}, as 
it corresponds to local equilibria for the Lorentz kinetic model (\ref{KntLrtz}).

As for the hydrodynamic (diffusion) limit of the periodic Lorentz gas with finite horizon, 
it was proved in \cite{BDG} that the case of an isotropic reflection law at the surface of 
each obstacle can be treated by PDE techniques, avoiding the heavy machinery from 
ergodic theory required to handle the case of specular reflection and developed by 
Bunimovich-Sinai --- see \cite{BuSi}, and also \cite{BuChSi}.

In the case of the Boltzmann-Grad limit however, changing the collision process does not affect
the result in Theorem \ref{TH-MN}, since the obstruction to using the linear Boltzmann equation
comes from particles that travel too far before encountering an obstacle for the first time. For
such particles, the nature of the collision process is obviously of no importance. Hence Theorem
\ref{TH-MN} holds verbatim if one replaces $\hat S^\eps_t(f^{in}\indc_{\Ups_\eps\times\bS^{D-1}})$
with the solution of
$$
\begin{aligned}
\d_tf_\eps+v\cdot\nabla_xf_\eps&=0\,,\quad x\in\Ups_\eps\,,\,\,v\in\bS^{D-1}\,,
\\
f_\eps(t)\rstr_{\Sigma_+^\eps}&=\cK(f_\eps(t)\rstr_{\Sigma_-^\eps})\,,
\\
f_\eps\rstr_{t=0}&=f^{in}\indc_{\Ups_\eps\times\bS^{D-1}}\,,
\end{aligned}
$$
where
$$
\Sigma_\pm^\eps=\{(x,v)\in\d\Ups_\eps\times\bS^{D-1}\,|\,\pm v\cdot n_x>0\}
$$
and $\cK$ is any linear operator from $L^\infty(\Sigma_-^\eps)$ to $L^\infty(\Sigma_+^\eps)$
that preserves the cone of positive functions.

\smallskip
By the same token, the same result as in Theorem \ref{TH-MN} holds without change if the 
obstacles are not assumed to be spherical, or even identical, but instead of arbitrary shapes, 
provided that they can be included in balls of radius $r=r_*\eps^{D/(D-1)}$ centered at the
points of $\eps\bZ^D$.

\bigskip
\noindent
{\bf Acknowledgements.} I express my gratitude to Profs. C. Bardos, C. Boldrighini and 
H.S. Dumas for helpful comments during the preparation of this paper.

\end{document}